\documentclass[11pt]{article}

\usepackage {amsmath, amssymb, graphics, emlines2}

\newtheorem{theorem}{Theorem}[section]
\newtheorem{lemma}[theorem]{Lemma}
\newtheorem{proposition}[theorem]{Proposition}

\begin{document}

\begin{center}\Large
\textbf{Tight contact structures on some bounded Seifert manifolds
with minimal convex boundary}
\end{center}
\begin{center}
\textbf{Fan Ding, Youlin Li and Qiang Zhang}\\
\end{center}

\hspace*{-0.2cm}\textbf{Abstract:} We classify positive tight
contact structures, up to isotopy fixing the boundary, on the
manifolds $N=M(D^{2}; r_1, r_2)$ with minimal convex boundary of
slope $s$ and Giroux torsion $0$ along $\partial N$, where
$r_1,r_2\in (0,1)\cap\mathbb{Q}$, in the following cases:

(1) $s\in(-\infty, 0)\cup[2, +\infty)$;

(2) $s\in[0, 1)$ and $r_1,r_2\in [1/2,1)$;

(3) $s\in[1, 2)$ and $r_1,r_2\in(0,\frac{1}{2})$;

(4) $s=\infty$ and $r_1=r_2=\frac{1}{2}$.

We also classify positive tight contact structures, up to isotopy
fixing the boundary, on $M(D^2;\frac{1}{2},\frac{1}{2})$ with
minimal convex boundary of arbitrary slope and Giroux torsion
greater than $0$ along the boundary.

\hspace*{-0.2cm}\textbf{Keywords:} \text{contact structure; bounded Seifert manifolds}

\hspace*{-0.2cm}\textbf{2010 Mathematical Subject Classification:}  \text{57M50; 53D10.}\\

\section{Introduction}

If $M$ is an oriented 3-manifold, a \emph{contact structure} on
$M$ is a completely non-integrable 2-plane distribution $\xi$
given as the kernel of a global 1-form $\alpha$ such that
$\alpha\wedge d\alpha\neq0$ at every point of $M$. Throughout this
paper, we assume the contact structures are \emph{positive}, i.e.,
given by one form $\alpha$ satisfying $\alpha\wedge d\alpha>0$,
and oriented.

Classification of tight contact structures on oriented 3-manifolds
is a fundamental problem in contact topology. See \cite{E},
\cite{Gi2}, \cite{Gi3}, \cite{Ho1} and \cite{Ho2}. The
classification of tight contact structures on small Seifert
manifolds has been the object of intensive study in the last few
years. See \cite{W}, \cite{GLS1}, \cite{GLS2} and \cite{G}. In
\cite{Ho1}, Honda classified tight contact structures on the solid
torus $S^1\times D^2$ and the thickened torus $T^{2}\times I$. In
\cite{T}, Tanya classified tight contact structures on
$\Sigma_{2}\times I$ where the boundary condition is specified by
a single, nontrivial separating dividing curve on each boundary
component. In this article, we classify tight contact structures
on some bounded Seifert manifolds.

Let $N$ be a small bounded Seifert manifold $M(D^{2}; r_1,r_2)$,
where $r_i\in (0,1)\cap\mathbb{Q},\, i=1,2$. We concentrate on
tight contact structures on $N$ with minimal convex boundary,
i.e., the number of dividing curves on $\partial N$ is $2$.
Suppose $s\in \mathbb{Q}$. Denote the greatest integer not greater
than $s$ by $[s]$. Let $s-[s]=\frac{b}{a}$, where $a>b\geq 0$ are
integers and $g.c.d. (a,b)=1$. If $\frac{1}{1-\frac{b}{a}}$ is not
an integer, then write
$\frac{1}{1-\frac{b}{a}}=a_1-\frac{1}{a_2-\cdots\frac{1}{a_{m-1}-\frac{1}{a_m}}}$,
where $a_j$'s are integers and $a_j\geq 2$ for $j\geq 1$. If
$\frac{1}{1-\frac{b}{a}}$ is an integer, then set
$a_1=\frac{1}{1-\frac{b}{a}}+1$, $a_2=1$ and $m=2$. Let
$r_3=\frac{1}{a_1-\frac{1}{a_2-\cdots\frac{1}{a_{m-1}-\frac{1}{a_m+1}}}}$.

\begin{theorem} \label{number} The number of tight contact structures on $N$ with minimal
convex boundary of slope $s$ and  Giroux torsion $0$ along
$\partial N$, up to isotopy fixing the boundary, is the number of
tight contact structures, up to isotopy, on the small Seifert
manifold $M(-1-[s];r_1,r_2,r_3)$ in the following cases:

\begin{enumerate}
\item $s\in (-\infty,0)$; \item $s\in [0,1)$ and $r_1,r_2\in
[\frac{1}{2},1)$; \item $s\in [1,2)$ and $r_1,r_2\in
(0,\frac{1}{2})$; \item $s\in [2,+\infty)$.

\end{enumerate}
\end{theorem}

The idea of the proof of Theorem 1.1 is roughly as follows.
Similar to the arguments in \cite{W}, \cite{GLS1}, \cite{GLS2}, we
get an upper bound for the number of tight contact structures, up
to isotopy fixing the boundary, on $N$ with given conditions. This
upper bound is the same as the number of tight contact structures,
up to isotopy, on $M(-1-[s];r_1,r_2,r_3)$. For any tight contact
structure $\eta$ on $M(-1-[s];r_1,r_2,r_3)$, we can decompose
$M(-1-[s];r_1,r_2,r_3)$ into $N$ and a solid torus $V_3$, and
isotope $\eta$ so that $\partial N$ is minimal convex with
dividing curves of slope $s$. When measured in the coordinates of
$\partial V_3$, this slope is $-1$. Thus, by the uniqueness of
tight contact structures, up to isotopy fixing the boundary, on a
solid torus with minimal convex boundary of slope $-1$, we
conclude that the number of tight contact structures, up to
isotopy, on $M(-1-[s];r_1,r_2,r_3)$ is less than or equal to the
number of tight contact structures, up to isotopy fixing the
boundary, on $N$ with given conditions.

Using the fact that a double cover of $M(D^2;
\frac{1}{2},\frac{1}{2})$ is the thickened torus $T^2\times I$ and
the classification of tight contact structures on $T^2\times I$,
we have the following classification.

\begin{theorem} \label{number1} (1) We can divide the set of tight
contact structures on $M(D^2; \frac{1}{2},\frac{1}{2})$ with
minimal convex boundary of slope $\infty$ and Giroux torsion $0$
along the boundary, up to isotopy fixing the boundary, into two
subsets. The tight contact structures in one subset are in 1-1
correspondence with $\mathbb{Z}$. The other subset contains two
elements.

(2) For any integer $t>0$ and any number $s\in \mathbb{Q}\cup\{
\infty\}$, there are exactly $2$ tight contact structures on
$M(D^2;\frac{1}{2},\frac{1}{2})$ with minimal convex boundary of
slope $s$ and Giroux torsion $t$ along the boundary, up to isotopy
fixing the boundary.
\end{theorem}

In Section 2, we give some preliminaries. In Section 3, we prove
Theorem \ref{number} in cases 1 and 2, and in Section 4, we prove
Theorem \ref{number} in cases 3 and 4. In Section 5, we prove
Theorem \ref{number1}. The reader is assumed to be familiar with
convex surfaces theory (cf. \cite{Gi1}, \cite{Ge}) and bypasses
(cf. \cite{Ho1}).

\section{Preliminaries}

For $r_1,r_2,r_3\in \mathbb{Q}\setminus \mathbb{Z}$, the Seifert
manifolds $M(D^2;r_1,r_2)$ and $M(r_1,r_2,r_3)$ are described as
follows. Let $\Sigma$ be an oriented pair of pants, and identify
each connected component of
$$-(\partial\Sigma\times S^{1})=T_{1}\cup T_{2}\cup T_{3}$$ with
$\mathbb{R}^{2}/\mathbb{Z}^{2}$, so that $\left( \begin{array}{c}
1 \\ 0 \end{array} \right)$ gives the direction of
$-\partial(\Sigma\times\{ 1\})$ and $\left(
\begin{array}{c} 0 \\ 1 \end{array} \right)$ gives the direction
of the $S^1$ factor. For $i=1,2,3$, let $V_{i}=D^{2}\times S^{1}$,
and identify $\partial V_i$  with
$\mathbb{R}^2/\mathbb{Z}^2$ so that $\left( \begin{array}{c} 1 \\
0 \end{array} \right)$ gives the direction of the meridian
$\partial (D^2\times \{ 1\})$ and $\left(
\begin{array}{c} 0 \\ 1 \end{array} \right)$ gives the direction
of the $S^1$ factor. Then $M(D^2;r_1,r_2)$ (respectively,
$M(r_1,r_2,r_3)$) is obtained from $\Sigma\times S^1$ by gluing
$V_i$ to $T_i$, $i=1,2$ (respectively, $i=1,2,3$), using the map
$\varphi_{i}:\partial V_{i}\to T_{i}$ defined by the matrix
$$\varphi_i=\left( \begin{array}{cc} p_i & u_i \\ -q_i & -v_i \end{array}
\right),$$ where $\frac{q_i}{p_i}=r_i$, $u_iq_i-p_iv_i=1$, and
$0<u_i<p_i$.

Note that if
$r_3=n+\frac{1}{a_1-\frac{1}{a_2-\cdots\frac{1}{a_{m-1}-\frac{1}{a_m+1}}}}$,
where $m\ge 2$, $n$ and $a_j$'s are integers, $a_j\ge 2$ for $1\le
j< m$ and $a_m\ge 1$, then $$\varphi_3=\left( \begin{array}{cc}
p_3 & u_3 \\ -q_3 & -v_3 \end{array} \right)
=\left( \begin{array}{cc} 1 & 0 \\
-n & 1
\end{array}\right) \left( \begin{array}{cc} a_1 & 1 \\ -1 & 0
\end{array} \right)\cdots \left( \begin{array}{cc} a_{m-1} & 1 \\ -1
& 0 \end{array} \right) \left( \begin{array}{cc} a_m+1 & 1 \\ -1 &
0 \end{array} \right)$$ and $$\left( \begin{array}{cc} a_1 & 1 \\
-1 & 0 \end{array} \right) \cdots \left( \begin{array}{cc} a_{m-1}
& 1
\\ -1 & 0
\end{array}\right) \left( \begin{array}{cc} a_m+1 &
1 \\ -1 & 0
\end{array}\right) \left( \begin{array}{c} -1 \\ 1
\end{array} \right) =\left( \begin{array}{c} -a\\ a-b \end{array}
\right) ,$$ where $a>b\ge 0$ are integers so that $g.c.d. (a,b)=1$
and
$\frac{b}{a}=1-\frac{1}{a_1-\frac{1}{a_2-\cdots\frac{1}{a_{m-1}-\frac{1}{a_m}}}}$.
Thus we have \begin{proposition} \label{slope-1}
In the notations above, $\varphi_3\left( \begin{array}{c} -1 \\
1
\end{array} \right) =\left( \begin{array}{c} -a \\ na+a-b
\end{array} \right)$. \hfill $\Box$
\end{proposition}

Let $n_1,n_2$ be integers. The Seifert manifolds $M(D^2;r_1,r_2)$
and $M(D^2;r_1+n_1,r_2+n_2)$ are orientation-preserving
diffeomorphic. This can be seen as follows. Let $f:\Sigma\times
S^1\to\Sigma\times S^1$ be an orientation-preserving diffeomorphism
such that $f$ sends each $T_i$ to itself and on each $T_i$, $f$ is
given by the matrix $f_i=\left(
\begin{array}{cc} 1 & 0 \\ -n_i & 1
\end{array} \right)$, where $n_3=-n_1-n_2$. ($f$ can be
constructed by using a smooth function $g:\Sigma\to SO(2)$ such
that for $x\in \Sigma$, $z\in S^1$, $f(x,z)=(x,g(x)z)$.) $f$ can
be extended to an orientation-preseving diffeomorphism, still
denoted by $f$, from $M(D^2;r_1,r_2)$ to $M(D^2;r_1+n_1,r_2+n_2)$.
Since $f_3=\left(
\begin{array}{cc} 1 & 0 \\ n_1+n_2 & 1
\end{array} \right)$, we have

\begin{proposition}\label{slopechange}
Under $f$, a simple closed curve of slope $s$ in $T_3$ of
$M(D^2;r_1,r_2)$ changes to a simple closed curve of slope
$s+n_1+n_2$ in $T_3$ of $M(D^2;r_1+n_1,r_2+n_2)$. \hfill $\Box$
\end{proposition}

Similarly, the Seifert manifolds $M(r_1,r_2,r_3)$ and
$M(r_1+n_1,r_2+n_2,r_3-n_1-n_2)$ are orientation-preserving
diffeomorphic. They are also denoted by
$M(e_0;r_1-[r_1],r_2-[r_2],r_3-[r_3])$, where
$e_0=[r_1]+[r_2]+[r_3]$.

On $T^{2}\times [0,1]\cong \mathbb{R}^{2}/\mathbb{Z}^{2}\times
[0,1]$ with coordinates $((x,y),t)$, consider
$\xi_{n}=\ker(\sin(\pi n t)dx+\cos(\pi n t)dy)$, with the boundary
adjusted so it becomes convex with two dividing curves on each
component, where $n\in \mathbb{Z}^{+}$. Let $(M,\xi)$ be a contact
3-manifold and $T\subset M$ an embedded torus. The \emph{Giroux
torsion} along $T$ is the supremum, over $n\in \mathbb{Z}^{+}$,
for which there exists a contact embedding $\phi: (T^{2}\times
[0,1], \xi_{n})\hookrightarrow (M,\xi)$, where
$\phi(T^{2}\times\{t\})$ is isotopic to $T$. (We set the Giroux
torsion to be $0$ if there is no such embedding). One can consult
\cite{HKM} for this definition.

The main invariant in the classification of tight contact
structures on Seifert manifolds is the maximal twisting number.
One can consult \cite{G} for the definition.

For the rest of the paper, $r_1,r_2\in (0,1)\cap\mathbb{Q}$,
$s\in\mathbb{Q}$, $a,b,a_j$ ($j=1,\ldots,m$) and $r_3$ are defined
as in the Introduction. For $i=1,2$, suppose
$-\frac{1}{r_i}=a_{0}^{i}-\frac{1}{a_{1}^{i}-\frac{1}{a_{2}^{i}-\cdots\frac{1}{a_{l_{i}-1}^{i}-\frac{1}{a_{l_{i}}^{i}}}}}$,
where $a_j^{i}$'s are integers and $a_{j}^{i}\leq-2$ for $j\geq0$.
When we consider the number of tight contact structures up to
isotopy or up to isotopy fixing the boundary, we usually omit the
phrase ``up to isotopy" or ``up to isotopy fixing the boundary".

\section{Proof of Theorem \ref{number} in cases 1 and 2}

For $i=1,2$, let $\varphi_{i}=\left(
                \begin{array}{cc}
                  p_{i} & u_{i} \\
                  -q_{i} & -v_{i} \\
                \end{array}
              \right)$, where $\frac{q_i}{p_i}=r_i$, $u_iq_i-p_iv_i=1$ and
              $0<u_i<p_i$.

Let $\xi$ be a tight contact structure on $N=M(D^2;r_1,r_2)$ with
minimal convex boundary of slope $s(T_3)$ and Giroux torsion $0$
along $\partial N$. We first isotope $\xi$ to make each $V_{i}$
($i=1,2$) a standard neighborhood of a Legendrian circle isotopic
to the $i$th singular fiber with twisting number $t_{i}<0$, i.e.,
$\partial V_{i}$ is convex with two dividing curves each of which
has slope $\frac{1}{t_{i}}$ when measured in the coordinates of
$\partial V_{i}$. Then, when measured in the coordinates of
$T_{i}$, the slope
$s_{i}=-\frac{q_{i}}{p_{i}}+\frac{1}{p_{i}(t_{i}p_{i}+u_{i})}<-\frac{q_{i}}{p_{i}}.$

The proof of the following lemma is similar to the proof of
\cite[Lemma 2.2]{W}.

\begin{lemma} \label{twistnumber0} On $M(D^{2}; r_1, r_2)$, if
$s(T_{3})\leq\max\{ r_1,r_2\}$, or if $0<s(T_{3})<1$ and
$r_i\geq\frac{1}{2}$ $(i=1,2)$, then any tight contact structure
with minimal convex boundary of slope $s(T_3)$ admits a vertical
Legendrian circle $L$ with twisting number $0$. \hfill $\Box$
\end{lemma}

Now suppose $s(T_3)<0$. Using the vertical Legendrian circle $L$,
we can thicken $V_{i}$ ($i=1,2$) to $V_{i}'$ such that $V_{i}'$'s
are pairwise disjoint, and $T_{i}'=\partial V_{i}'$ is a minimal
convex torus with vertical dividing curves when measured in
coordinates of $T_{i}$. Also, we can thicken $T_{3}$ to
$L_{3}=T_{3}\times[0,1]$ such that $T_{3}\times\{0\}=T_{3}$ and
$T_{3}\times\{1\}=T_{3}'$ is a minimal convex torus with vertical
dividing curves when measured in the coordinates of $T_{3}$.
Choose $t_i\ll -1$ so that $-\infty<\frac{1}{a_{0}^{i}+1}<s_{i}$
for $i=1,2$. By \cite[Proposition 4.16]{Ho1}, for $i=1,2$, there
exists a minimal convex torus $T_{i}''$ in the interior of
$V_{i}'\setminus V_{i}$ isotopic to $T_{i}$ that has dividing
curves of slope $\frac{1}{a_{0}^{i}+1}$. Let $V_{i}''$ be the
solid torus bounded by $T_{i}''$, and $\Sigma''\times
S^{1}=N\setminus(V_{1}''\cup V_{2}'')$.

First we consider $V_{1}''$ and $V_{2}''$. Since
$\varphi_{i}^{-1}\left(
                                        \begin{array}{c}
                                          a_{0}^{i}+1 \\
                                          1 \\
                                        \end{array}
                                      \right)=\left(
                                                \begin{array}{c}
                                                  -(a_{0}^{i}+1)v_{i}-u_{i} \\
                                                  (a_{0}^{i}+1)q_{i}+p_{i} \\
                                                \end{array}
                                              \right)$ (here $u_i,v_i$ correspond
respectively to $-u_i,-v_i$ in the proof of \cite[Theorem
1.6]{W}), the dividing curves of $T_{i}''$ ($i=1,2$) have slope
$-\frac{(a_{0}^{i}+1)q_{i}+p_{i}}{(a_{0}^{i}+1)v_{i}+u_{i}}$ when
measured in the coordinates of $\partial V_{i}$. By a similar
argument as in the proof of \cite[Theorem 1.6]{W}, there are
exactly $\prod_{j=1}^{l_{i}}|a_{j}^{i}+1|$ tight contact
structures on $V_{i}''$ that satisfy the given boundary condition.

Then we consider $N\setminus(V_{1}''\cup V_{2}'')=\Sigma''\times
S^{1}$. Let $L_{i}$ ($i=1,2$) be the thickened torus which is
bounded by $T_{i}'$ and $T_{i}''$, then $L_{i}$ has boundary
slopes $\infty$ and $\frac{1}{a_{0}^{i}+1}$. By \cite[Theorem
2.2]{Ho1}, there are exactly $|a_{0}^{i}|$ minimally twisting
tight contact structures on $L_{i}$ that satisfy the given
boundary condition. The two boundary slopes of the thickened
torus $L_{3}$ are $\infty$ and $s(T_3)$ respectively.\\

\textbf{Case 1.} $s\in (-\infty, 0)$.

We divide it into two subcases.

\textbf{Case 1(a).} $s\in (-\infty,-1).$

Let $s(T_3)=s$. We decompose $L_3$ into $m$ continued fraction
blocks (some blocks may be invariant neighborhoods of convex
tori). The first continued fraction block has two boundary slopes
$\infty$ and $[s]+1-\frac{1}{a_{1}-1}$, the second continued
fraction block has two boundary slopes $[s]+1-\frac{1}{a_{1}-1}$
and $[s]+1-\frac{1}{a_{1}-\frac{1}{a_{2}-1}}$, $\ldots$, the $m$th
continued fraction block has two boundary slopes
$[s]+1-\frac{1}{a_{1}-\frac{1}{a_{2}-\cdots\frac{1}{a_{m-1}-1}}}$
and
$s(T_{3})=s=[s]+1-\frac{1}{a_{1}-\frac{1}{a_{2}-\cdots\frac{1}{a_{m-1}-\frac{1}{a_{m}}}}}$.
By shuffling, there are at most
$a_{1}(a_{2}-1)\ldots(a_{m-1}-1)a_{m}$ minimally twisting tight
contact structures on $L_{3}$.  By a similar argument as in the
first paragraph of \cite[page 241]{W}, the upper bound of the
number of tight contact structures on $N$ with minimal convex
boundary of slope $s$ and Giroux torsion $0$ along $\partial N$ is
$|a_{0}^{1}
a_{0}^{2}\prod_{i=1}^{2}\prod_{j=1}^{l_{i}}(a_{j}^{i}+1)|a_{1}(a_{2}-1)\ldots(a_{m-1}-1)a_{m}$.

Consider the closed Seifert manifold $M(r_1,r_2,-1-[s]+r_3)$.
Since $-1-[s]>0$, by \cite[Theorem 1.6]{W}, it admits $|a_{0}^{1}
a_{0}^{2}\prod_{i=1}^{2}\prod_{j=1}^{l_{i}}(a_{j}^{i}+1)|a_{1}(a_{2}-1)\ldots(a_{m-1}-1)a_{m}$
tight contact structures. By \cite[Theorem 1.3]{W}, for any tight
contact structure $\eta$ on $M(r_1,r_2, -1-[s]+r_3)$, there is a
vertical Legendrian circle with twisting number $0$. We isotope
$\eta$ so that there is a vertical Legendrian circle $L$ with
twist number $0$ in the interior of $\Sigma\times S^1$, and
$V_{3}$ is a standard neighborhood of a Legendrian circle isotopic
to the $3$rd singular fiber with twisting number $t<0$, i.e.,
$\partial V_{3}$ is convex with two dividing curves each of which
has slope $\frac{1}{t}$ when measured in the coordinates of
$\partial V_{3}$. Let $\varphi_3=\left( \begin{array}{cc} p_3 &
u_3 \\ -q_3 & -v_3
\end{array} \right)$, where $\frac{q_3}{p_3}=-1-[s]+r_3$, $u_3q_3-p_3v_3=1$ and
$0<u_3<p_3$. Then, when measured in the coordinates of $T_{3}$,
the slope $s_3=-\frac{q_3}{p_3}+\frac{1}{p_3(t_3p_3+u_3)}$. Using
$L$, we can thicken $V_{3}$ to $V_{3}'$, such that
$T_{3}'=\partial V_{3}'$ is a minimal convex torus with vertical
dividing curves when measured in the coordinates of $T_{3}$. Since
$1-\frac{b}{a}=\frac{1}{a_{1}-\frac{1}{a_{2}-\cdots\frac{1}{a_{m-1}-\frac{1}{a_{m}}}}}>$
$\frac{1}{a_{1}-\frac{1}{a_{2}-\cdots\frac{1}{a_{m-1}-\frac{1}{a_{m}+1}}}}=r_3$,
we have $-\infty<s=[s]+\frac{b}{a}<[s]+1-r_3$. Thus $-\infty
<s<s_3$ for sufficiently small $t$. By \cite[Proposition
4.16]{Ho1}, there exists a minimal convex torus $T_{3}''$ in the
interior of $V_{3}'\setminus V_{3}$ isotopic to $T_{3}$ that has
dividing curves of slope $s$. Thus we can isotopy $\eta$ so that
$T_3$ is minimal convex with dividing curves of slope $s$ when
measured in the coordinates of $T_3$. Note that
$M(r_1,r_2,-1-[s]+r_3)$ has a decomposition $N\cup_{\varphi_3}
V_{3}$. By Proposition \ref{slope-1}, $\varphi_{3}^{-1}\left(
                    \begin{array}{c}
                      -a \\
                      -[s]a-b \\
                    \end{array}
                  \right)=
                  \left(
                            \begin{array}{c}
                              -1 \\
                            1 \\
                            \end{array}
                            \right)$.
Thus the slope of the dividing curves on $\partial V_{3}$ is $-1$
when measured in the coordinates of $\partial V_{3}$. There is
exactly one tight contact structure on $V_{3}$ with minimal convex
boundary of slope $-1$. Note that $\eta$, when restricted to $N$,
has Giroux torsion $0$ along $\partial N$. Hence the number of
tight contact structures on $N$ with given conditions is at least
$|a_{0}^{1}
a_{0}^{2}\prod_{i=1}^{2}\prod_{j=1}^{l_{i}}(a_{j}^{i}+1)|a_{1}(a_{2}-1)\ldots(a_{m-1}-1)a_{m}$.

Therefore, there are exactly $|a_{0}^{1}
a_{0}^{2}\prod_{i=1}^{2}\prod_{j=1}^{l_{i}}(a_{j}^{i}+1)|a_{1}(a_{2}-1)\ldots(a_{m-1}-1)a_{m}$
tight contact structures on $N$ with minimal convex boundary of
slope $s$ and Giroux torsion $0$ along $\partial N$.

\textbf{Case 1(b).} $s\in [-1,0)$.

Let $s(T_3)=s$. Note that the outermost continued fraction block
of $L_{3}$ has two boundary slopes $\infty$ and
$-\frac{1}{a_{1}-1}$, and hence contains $a_{1}-1$ basic slices.
By a similar argument as in the proof of \cite[Theorem 2.4]{GLS1},
there are at most $[a_{0}^{1} a_{0}^{2}a_{1}-(a_{0}^{1}+1)
(a_{0}^{2}+1)(a_{1}-1)](a_{2}-1)\ldots(a_{m-1}-1)a_{m}\prod_{i=1}^{2}\prod_{j=1}^{l_{i}}|a_{j}^{i}+1|$
tight contact structures on $N$ with minimal convex boundary of
slope $s$ and Giroux torsion $0$ along $\partial N$.

Consider the small Seifert manifold $M(r_1,r_2,r_3)$. By
\cite[Theorem 1.1]{GLS1}, it admits exactly $[a_{0}^{1}
a_{0}^{2}a_{1}-(a_{0}^{1}+1)
(a_{0}^{2}+1)(a_{1}-1)](a_{2}-1)\ldots(a_{m-1}-1)a_{m}\prod_{i=1}^{2}\prod_{j=1}^{l_{i}}|a_{j}^{i}+1|$
tight contact structures. Let $\varphi_3=\left( \begin{array}{cc}
p_3 & u_3 \\ -q_3 & -v_3
\end{array} \right)$, where $\frac{q_3}{p_3}=r_3$, $u_3q_3-p_3v_3=1$ and
$0<u_3<p_3$. Note that $M(r_1,r_2,r_3)$ has a decomposition
$N\cup_{\varphi_3} V_{3}$. By a similar argument  as in Case 1(a),
for any tight contact structure $\eta$ on $M(r_1,r_2,r_3)$, we can
isotopy $\eta$ so that $T_3$ is minimal convex with dividing
curves of slope $s$ when measured in the coordinates of $T_3$. By
Proposition \ref{slope-1}, $\varphi_{3}^{-1}\left(
                    \begin{array}{c}
                      -a \\
                      a-b \\
                    \end{array}
                  \right)=
                  \left(
                            \begin{array}{c}
                              -1 \\
                            1 \\
                            \end{array}
                            \right)$.
Thus the slope of the dividing curves on $\partial V_{3}$ is $-1$
when measured in the coordinates of $\partial V_{3}$. Similar to
Case 1(a),  we conclude that the number of tight contact
structures on $N$ with given conditions is at least $[a_{0}^{1}
a_{0}^{2}a_{1}-(a_{0}^{1}+1)
(a_{0}^{2}+1)(a_{1}-1)](a_{2}-1)\ldots(a_{m-1}-1)a_{m}$
$\prod_{i=1}^{2}\prod_{j=1}^{l_{i}}|a_{j}^{i}+1|$.

Therefore, there are exactly $[a_{0}^{1}
a_{0}^{2}a_{1}-(a_{0}^{1}+1)
(a_{0}^{2}+1)(a_{1}-1)](a_{2}-1)\ldots(a_{m-1}-1)a_{m}$
$\prod_{i=1}^{2}\prod_{j=1}^{l_{i}}|a_{j}^{i}+1|$ tight contact
structures on $N$ with the given boundary condition and Giroux
torsion $0$ along $\partial N$.

\textbf{Case 2.} $s\in [0,1)$ and $r_1,r_2\in [\frac{1}{2},1)$.

By Lemma \ref{twistnumber0}, any tight contact structure on
$M(D^2;r_1,r_2)$ with minimal convex boundary of slope $s$
contains a Legendrian vertical circle with twisting number $0$.

By Proposition \ref{slopechange}, a tight contact structure on
$M(D^2;r_1,r_2)$ with minimal convex boundary of slope $s$
corresponds to a tight contact structure on $M(D^2;-1+r_1,r_2)$
with minimal convex boundary of slope $s-1$. We consider tight
contact structures on $M(D^2;-1+r_1,r_2)$ with minimal convex
boundary of slope $s-1$. Without loss of generality, assume that
$r_1\ge r_2$.

Suppose $\xi_0$ is a tight contact structure on $N_0=M(D^{2};
-1+r_1, r_2)$ with minimal convex boundary of slope $s-1$ and
Giroux torsion $0$ along $\partial N_0$. Using a vertical
Legendrian circle with twisting number $0$, we can thicken
standard neighborhoods of two Legendrian singular fibers to
$U_{1}$ and $U_{2}$ such that the slopes of the dividing curves on
$\partial U_{1}$ and $\partial U_{2}$ are $\infty$ when measured
in the coordinates of $T_1$ and $T_2$, respectively. We can
thicken $T_{3}$ to a thickened torus $L_{3}$ so that the slope of
the other boundary component of $L_{3}$ is $\infty$.

Consider the closed Seifert manifold $M(-1+r_1,r_2,r_3)$. Let
$\varphi_3=\left( \begin{array}{cc} p_3 & u_3 \\ -q_3 & -v_3
\end{array} \right)$, where $\frac{q_3}{p_3}=r_3$, $u_3q_3-p_3v_3=1$ and
$0<u_3<p_3$. $M(r_1,r_2,r_3)$ has a decomposition
$N\cup_{\varphi_3} V_{3}$. By Proposition \ref{slope-1},
$\varphi_{3}^{-1}\left(
                    \begin{array}{c}
                      -a \\
                      a-b \\
                    \end{array}
                  \right)=
                  \left(
                            \begin{array}{c}
                              -1 \\
                            1 \\
                            \end{array}
                            \right)$.
Thus we can find layers $N_j^i$ of $M(-1+r_1,r_2,r_3)$ in
\cite[page 1432]{GLS2} in $(N_0,\xi_0)$. Note that $L_3$
corresponds to $U_3\setminus V_3'$ in \cite[page 1432]{GLS2}. We
can obtain similar results as in \cite[Propositions 6.1 and
6.3]{GLS2} for $(N_0,\xi_0)$. Since $\frac{1}{r_3}$ is not an
integer, $L_3$ contains at least two layers and we can obtain a
similar result as in \cite[Proposition 6.4]{GLS2} for
$(N_0,\xi_0)$ (we only encounter Case 2 in the proof of
\cite[Proposition 6.4]{GLS2}). Therefore we obtain an upper bound
of the number of tight contact structures on $N_0$ with given
conditions, which is the same as the number of tight contact
structures on $M(-1+r_1,r_2,r_3)$.

By \cite[Proposition 5.1]{GLS2}, for any tight contact structure
$\eta$ on $M(-1+r_1, r_2,r_3)$, there is a Legendrian vertical
circle with twisting number $0$. Then similar to Case 1(a), we can
isotopy $\eta$ so that $T_3=\partial V_3$ is minimal convex with
dividing curves of slope $-1$ when measured in the coordinates of
$\partial V_3$. Then we conclude that the number of  tight contact
structures on $M(-1+r_1,r_2,r_3)$ is less than or equal to the
number of tight contact structures on $N_0$ with minimal convex
boundary of slope $s-1$ and Giroux torsion $0$ along $\partial
N_0$.

Therefore, the number of  tight contact structures on $N$ with
given conditions is exactly the number of tight contact structures
on $M(-1;r_1,r_2,r_3)$.

\section{Proof of Theorem \ref{number} in cases 3 and 4}

By Proposition \ref{slopechange}, a tight contact structure on
$M(D^2;r_1,r_2)$ with minimal convex boundary of slope $s$
corresponds to a tight contact structure on $M(D^2;r_1-1,r_2-1)$
with minimal convex boundary of slope $s-2$.

Now consider the manifold $M=M(D^2;r_1-1,r_2-1)$ and let
$s(T_3)=s-2$. For $i=1,2$,
$r_i-1=-1-\frac{1}{a_{0}^{i}-\frac{1}{a_{1}^{i}-\cdots\frac{1}{a_{l_{i}-1}^{i}-\frac{1}{a_{l_{i}}^{i}}}}}$.
Suppose $r_i-1=-\frac{q_i}{p_i}$, where $p_i,q_i$ are integers,
$0<q_i<p_i$ and $g.c.d. (p_i,q_i)=1$. Let $\varphi_{i}=\left(
                \begin{array}{cc}
                  p_{i} & u_{i} \\
                  q_{i} & v_{i} \\
                \end{array}
              \right)$, where $p_{i}> u_{i}>0$ and $p_{i}v_{i}-q_{i}u_{i}=1$.

Let $\xi$ be a tight contact structure on $M$ with minimal convex
boundary of slope $s(T_3)$ and Giroux torsion $0$ along $\partial
M$. The proof of the following lemma is similar to the proof of
\cite[Theorem 1.4]{W}.

\begin{lemma} \label{notwistnumber0}
On $M=M(D^{2}; -\frac{q_{1}}{p_{1}}, -\frac{q_{2}}{p_{2}})$, if
$s(T_{3})\geq -1$, then any tight contact structure with minimal
convex boundary of slope $s(T_3)$ and Giroux torsion $0$ along
$\partial M$ does not admit Legendrian vertical circles with
twisting number $0$. \hfill $\Box$
\end{lemma}

The proof of the following lemma is similar to the proof of the
corresponding result contained in the proof of \cite[Theorem
1.7]{W}.

\begin{lemma} \label{twistnumber-1}
On $M=M(D^{2}; -\frac{q_{1}}{p_{1}}, -\frac{q_{2}}{p_{2}})$, if
$s(T_{3})\geq 0$, or if $-1\leq s(T_{3})<0$ and
$\frac{q_{i}}{p_{i}}>\frac{1}{2}$ $(i=1,2)$, then the maximal
twisting number of Legendrian vertical circles in $(M,\xi)$ is
$-1$.
\end{lemma}

Now assume that $s(T_3)\ge 0$, or $-1\leq s(T_{3})<0$ and
$\frac{q_{i}}{p_{i}}>\frac{1}{2}$ $(i=1,2)$. By Lemma
\ref{twistnumber-1}, after an isotopy of $\xi$, we can find a
Legendrian vertical circle $L$ in the interior of $\Sigma\times
S^{1}$ with twisting number $-1$. Then we make each $V_{i}$
($i=1,2$) a standard neighborhood of a Legendrian circle which is
isotopic to the $i$th singular fiber with twisting number
$t_{i}<-2$, i.e., $\partial V_{i}$ is convex with two dividing
curves each of which has slope $\frac{1}{t_{i}}$ when measured in
the coordinates of $\partial V_{i}$. Let $s_{i}$ be the slope of
the dividing curves of $T_{i}=\varphi_{i}(\partial V_{i})$
measured in the coordinates of $T_{i}$. Then we have
$s_{i}=\frac{q_{i}t_{i}+v_{i}}{p_{i}t_{i}+u_{i}}=\frac{q_{i}}{p_{i}}+\frac{1}{p_{i}(p_{i}t_{i}+u_{i})}$.
The fact that $t_{i}<-2$ implies that
$0<s_{i}<\frac{q_{i}}{p_{i}}$. In particular, if
$\frac{q_{i}}{p_{i}}>\frac{1}{2}$ $(i=1,2)$, then $\frac{1}{2}<
s_{i}<\frac{q_{i}}{p_{i}}$.

We can assume that $T_{i}=\varphi_{i}(\partial V_{i})$ ($i=1,2$)
and $T_{3}$ have Legendrian rulings of slope $\infty$ when
measured in the coordinates of $T_i$ and $T_3$, respectively.
Using $L$, we can thicken $V_i$ to $V_i'$ ($i=1,2$) and $T_3$ to
$L_3$ to get a decomposition, $M=(\Sigma'\times
S^{1})\cup(V_{1}'\cup V_{2}'\cup L_{3})$, such that
$T_{i}'=\partial V_{i}'$ $(i=1,2)$ has two dividing curves of
slope $[s_{i}]=0$ when measured in the coordinates of $T_{i}$, and
the thickened torus $L_{3}$ has two boundary slopes $[s(T_3)]$ and
$s(T_3)$ (cf. the proof of \cite[Theorem 1.7]{W}).

By \cite[Lemma 5.1]{Ho2}, there are exactly $2+[s(T_3)]$ tight
contact structures on $\Sigma'\times S^{1}$ satisfying the
boundary condition and admitting no Legendrian vertical circles
with twisting number $0$.

The slope of the dividing curves on $\partial V_{i}'$ is
$-\frac{q_{i}}{v_{i}}$ when measured in the coordinates of
$\partial V_{i}$. So by a similar argument as in the proof of
\cite[Theorem 1.7]{W}, there are exactly
$\prod_{j=0}^{l_{i}}|a_{j}^{i}+1|$ tight contact structures on
$V_{i}'$ satisfying such boundary condition.

We consider Case 4 first.

\textbf{Case 4.} $s\in [2,+\infty)$.

We decompose $L_3$ into $m$ continued fraction blocks (some blocks
may be invariant neighborhoods of convex tori). The first one has
boundary slopes $[s]-2$ and $[s]-1-\frac{1}{a_{1}-1}$, the second
one has boundary slopes $[s]-1-\frac{1}{a_{1}-1}$ and
$[s]-1-\frac{1}{a_{1}-\frac{1}{a_{2}-1}}$, $\ldots$, the last one
has boundary slopes
$[s]-1-\frac{1}{a_{1}-\frac{1}{a_{2}-\cdots\frac{1}{a_{m-1}-1}}}$
and
$[s]-1-\frac{1}{a_{1}-\frac{1}{a_{2}-\cdots\frac{1}{a_{m-1}-\frac{1}{a_{m}}}}}=s(T_3)$.
By shuffling, there are at most
$(a_{1}-1)(a_{2}-1)\ldots(a_{m-1}-1)a_{m}$ minimally twisting
tight contact structures on $L_{3}$ (except $s=[s]$). So there are
at most
$[s]\prod_{i=1}^{2}\prod_{j=0}^{l_{i}}|a_{j}^{i}+1|(a_{1}-1)(a_{2}-1)\ldots(a_{m-1}-1)a_{m}$
tight contact structures on $M$ with minimal convex boundary of
slope $s-2$ and Giroux torsion $0$ along $\partial M$.

Consider the small closed Seifert manifold $M(r_1-1,r_2
-1,-[s]+1+r_3)=M(-\frac{q_{1}}{p_{1}}, -\frac{q_{2}}{p_{2}},
-[s]+1+r_3)$. Since $[s]\ge 2$, applying \cite[Theorem 1.7]{W},
there are exactly
$[s]\prod_{i=1}^{2}\prod_{j=0}^{l_{i}}|a_{j}^{i}+1|(a_{1}-1)(a_{2}-1)\ldots(a_{m-1}-1)a_{m}$
tight contact structures on $M(-\frac{q_{1}}{p_{1}},
-\frac{q_{2}}{p_{2}}, -[s]+1+r_3)$. According to the proof of
\cite[Theorem 1.7]{W}, for any tight contact structure $\eta$ on
$M(-\frac{q_{1}}{p_{1}}, -\frac{q_{2}}{p_{2}}, -[s]+1+r_3)$, the
maximal twisting number of a Legendrian vertical circle is $-1$.
After an isotopy of $\eta$,  we can find a vertical Legendrian
circle $L'$ with twist number $-1$ in the interior of
$\Sigma\times S^1$ and make  $V_{3}$ a standard neighborhood of a
Legendrian circle isotopic to the $3$rd singular fiber with
twisting number $t<0$, i.e., $\partial V_{3}$ is convex with two
dividing curves each of which has slope $\frac{1}{t}$ when
measured in the coordinates of $\partial V_{3}$. Let
$\varphi_3=\left(
\begin{array}{cc} p_3 & u_3 \\ q_3 & v_3
\end{array} \right)$, where $\frac{q_3}{p_3}=[s]-1-r_3$, $p_3v_3-u_3q_3=1$ and
$0<u_3<p_3$. Then, when measured in the coordinates of $T_{3}$,
the slope $s_{3}=\frac{q_3}{p_3}+\frac{1}{p_3(t_3p_3+u_3)}$.

Using $L'$, we can thicken $V_{3}$ to $V_{3}'$ such that
$T_{3}'=\partial V_{3}'$ has two dividing curves of slope $[s]-2$.
Since $[s]-2\le s-2<[s]-1-r_3$, $[s]-2\le s-2<s_3$ for
sufficiently small $t$. By \cite[Proposition 4.16]{Ho1}, there is
a convex torus $T_{3}''$ in the interior of $V_3'\setminus V_3$
which is parallel to $T_{3}$ and has two dividing curves of slope
$s(T_3)=s-2$. Thus we can isotopy $\eta$ so that $T_3$ is minimal
convex with dividing curves of slope $s(T_3)$ when measured in the
coordinates of $T_3$. $M(r_1-1,r_2-1,-[s]+1+r_3)$ has a
decomposition $M\cup_{\varphi_3} V_{3}$. By Proposition
\ref{slope-1}, $\varphi_{3}^{-1}\left(
                    \begin{array}{c}
                      -a \\
                      -([s]-2)a-b \\
                    \end{array}
                  \right)=
                  \left(
                            \begin{array}{c}
                              -1 \\
                            1 \\
                            \end{array}
                            \right)$.
Thus the slope of the dividing curves on $\partial V_{3}$ is $-1$
when measured in the coordinates of $\partial V_{3}$. Similar to
the argument in Case 1(a), the number of tight contact structures
on $M$ with given conditions is at least
$[s]\prod_{i=1}^{2}\prod_{j=0}^{l_{i}}|a_{j}^{i}+1|(a_{1}-1)(a_{2}-1)\ldots(a_{m-1}-1)a_{m}$.

Therefore, there are exactly
$[s]\prod_{i=1}^{2}\prod_{j=0}^{l_{i}}|a_{j}^{i}+1|(a_{1}-1)(a_{2}-1)\ldots(a_{m-1}-1)a_{m}$
tight contact structures on $M$ with minimal convex boundary of
slope $s-2$ and Giroux torsion $0$ along $\partial M$.

\textbf{Case 3.} $s\in [1,2)$ and $r_1,r_2\in (0,\frac{1}{2})$.

Since $r_i\in (0,\frac{1}{2})$ ($i=1,2$), $\frac{q_i}{p_i}=1-r_i >
\frac{1}{2}$ ($i=1,2$). Similar to Case 4, there are at most
$\prod_{i=1}^{2}\prod_{j=0}^{l_{i}}|a_{j}^{i}+1|(a_{1}-1)(a_{2}-1)\ldots(a_{m-1}-1)a_{m}$
tight contact structures on $M$ with given conditions.

Consider the small closed Seifert manifold $M(r_1-1,r_2-1,r_3)
=M(-\frac{q_1}{p_1},-\frac{q_2}{p_2},r_3)$. We claim that this
small closed Seifert manifold is an $L$-space (see \cite{LS} for
the definition). Note that since $r_1+r_2+r_3-2\neq 0$,
$M(r_1-1,r_2-1,r_3)$ is a rational homology sphere. By
\cite[Theorem 1.1]{LS}, it suffices to show that
$-M(-\frac{q_{1}}{p_{1}}, -\frac{q_{2}}{p_{2}}, r_3)=M(-1;
\frac{q_{1}}{p_{1}}, \frac{q_{2}}{p_{2}},1-r_3)$ carries no
positive, transverse contact structures. Suppose otherwise, by
\cite[Theorem 4.5]{G}, there are integers $h_{1}$, $h_{2}$,
$h_{3}$ and $k>0$, such that (1)
$\frac{h_{1}}{k}<-\frac{q_{1}}{p_{1}}$,
$\frac{h_{2}}{k}<-\frac{q_{2}}{p_{2}}$, $\frac{h_{3}}{k}<r_3-1$,
and (2) $\frac{h_{1}+h_{2}+h_{3}}{k}=-1-\frac{1}{k}$. Let $n$ be a
positive integer such that $1-r_3\geq\frac{1}{n}$. Combining (1)
and (2), we have
$-1-\frac{1}{k}<-\frac{q_{1}}{p_{1}}-\frac{q_{2}}{p_{2}}+r_3-1<
-1-\frac{1}{n}$. So $1 \leq k\leq n-1$. If $k$ is even, then
$h_{1}\leq\frac{-k}{2}-1$, $h_{2}\leq\frac{-k}{2}-1$ and
$h_{3}\leq -1$. Thus  we have
$-1-\frac{1}{k}=\frac{h_{1}+h_{2}+h_{3}}{k}\leq
\frac{-k-3}{k}=-1-\frac{3}{k}$. This is absurd. If $k$ is odd,
then $h_{1}\leq[\frac{-k}{2}]=\frac{-k-1}{2}$,
$h_{2}\leq[\frac{-k}{2}]=\frac{-k-1}{2}$ and $h_{3}\leq -1$. Thus
we have $-1-\frac{1}{k}=\frac{h_{1}+h_{2}+h_{3}}{k}\leq
\frac{2(\frac{-k-1}{2})-1}{k}=-1-\frac{2}{k}$. This is also
absurd.

By \cite[Theorem 1.3]{G}, there are exactly
$\prod_{i=1}^{2}\prod_{j=0}^{l_{i}}|a_{j}^{i}+1|(a_{1}-1)(a_{2}-1)\ldots(a_{m-1}-1)a_{m}$
tight contact structures on $M(r_1-1,r_2-1, r_3)$. Moreover, in
each of these tight contact structures, there is a Legendrian
vertical circle with twisting number $-1$. By a similar argument
as in Case 4, there are at least
$\prod_{i=1}^{2}\prod_{j=0}^{l_{i}}|a_{j}^{i}+1|(a_{1}-1)(a_{2}-1)\ldots(a_{m-1}-1)a_{m}$
tight contact structures on $M$ with given conditions. So, in this
case, there are exactly
$\prod_{i=1}^{2}\prod_{j=0}^{l_{i}}|a_{j}^{i}+1|(a_{1}-1)(a_{2}-1)\ldots(a_{m-1}-1)a_{m}$
tight contact structures on $M$ with minimal convex boundary of
slope $s-2$ and Giroux torsion $0$ along $\partial M$.

\section{Proof of Theorem \ref{number1}}

Consider the thickened torus $S^{1}\times S^{1}\times I$, where
$I=[0,1]$. Denote by $(x,y,z)$ the coordinates of $S^{1}\times
S^{1}\times I$, where $x\in \mathbb{R}/2\pi \mathbb{Z}$, $y\in
\mathbb{R}/2\pi \mathbb{Z}$ and $z\in[0,1]$. For convenience, let
$x\in[-\pi,\pi]$, $y\in[-\pi,\pi]$, and $-\pi$ is identified with
$\pi$. Let $N$ be the quotient space of $S^1\times S^1\times I$ by
identifying $(x,y,z)$ with $(x+\pi,-y,1-z)$. Let $p:S^{1}\times
S^{1}\times I\longrightarrow N$ be the covering projection. The
covering transformation $(x,y,z)\mapsto (x+\pi,-y,1-z)$ is denoted
by $\tau$. $N$ can be identified with
$M(D^2;-\frac{1}{2},\frac{1}{2})$ (the images under $p$ of
$S^1\times \{ 0\}\times \{\frac{1}{2}\}$ and $S^1\times
\{\pi\}\times\{\frac{1}{2}\}$ correspond to the singular fibers).
On the boundary $T_3=-\partial N=S^1\times S^1\times\{0\}$,
$S^1\times\{ pt\}\times\{ 0\}$ gives the fiber direction (i.e.,
corresponds to $\left(
\begin{array}{c} 0 \\ 1 \end{array} \right)$ in the notation of
Section 2) and $\{ pt\}\times S^1 \times\{0\}$ corresponds to
$\left( \begin{array}{c} -1 \\
0 \end{array} \right)$ in the notation of Section 2. By
Proposition \ref{slopechange}, a simple closed curve of slope $s$
in $T_3$ of $N=M(D^2;-\frac{1}{2},\frac{1}{2})$ corresponds to a
simple closed curve of slope $s+1$ in $T_3$ of
$M(D^2;\frac{1}{2},\frac{1}{2})$.

We regard $N$ as the quotient space of a thickened cylinder $[0,
\pi]\times S^{1} \times  I$ by identifying $(0, y, z)$ with $(\pi,
-y, 1-z)$.  See Figure 1. The coordinates of the four points $P$,
$Q$, $R$ and $S$ at the left end are $(0,-\frac{\pi}{2}, 0)$,
$(0,-\frac{\pi}{2}, 1)$, $(0, \frac{\pi}{2}, 1)$ and $(0,
\frac{\pi}{2}, 0)$ respectively. The coordinates of the two points
$X$ and $Y$ at the left end are $(0, 0, \frac{1}{2})$ and  $(0,
\pi, \frac{1}{2})$ respectively.

\begin{center}
\scalebox{0.25}[0.25]{\includegraphics {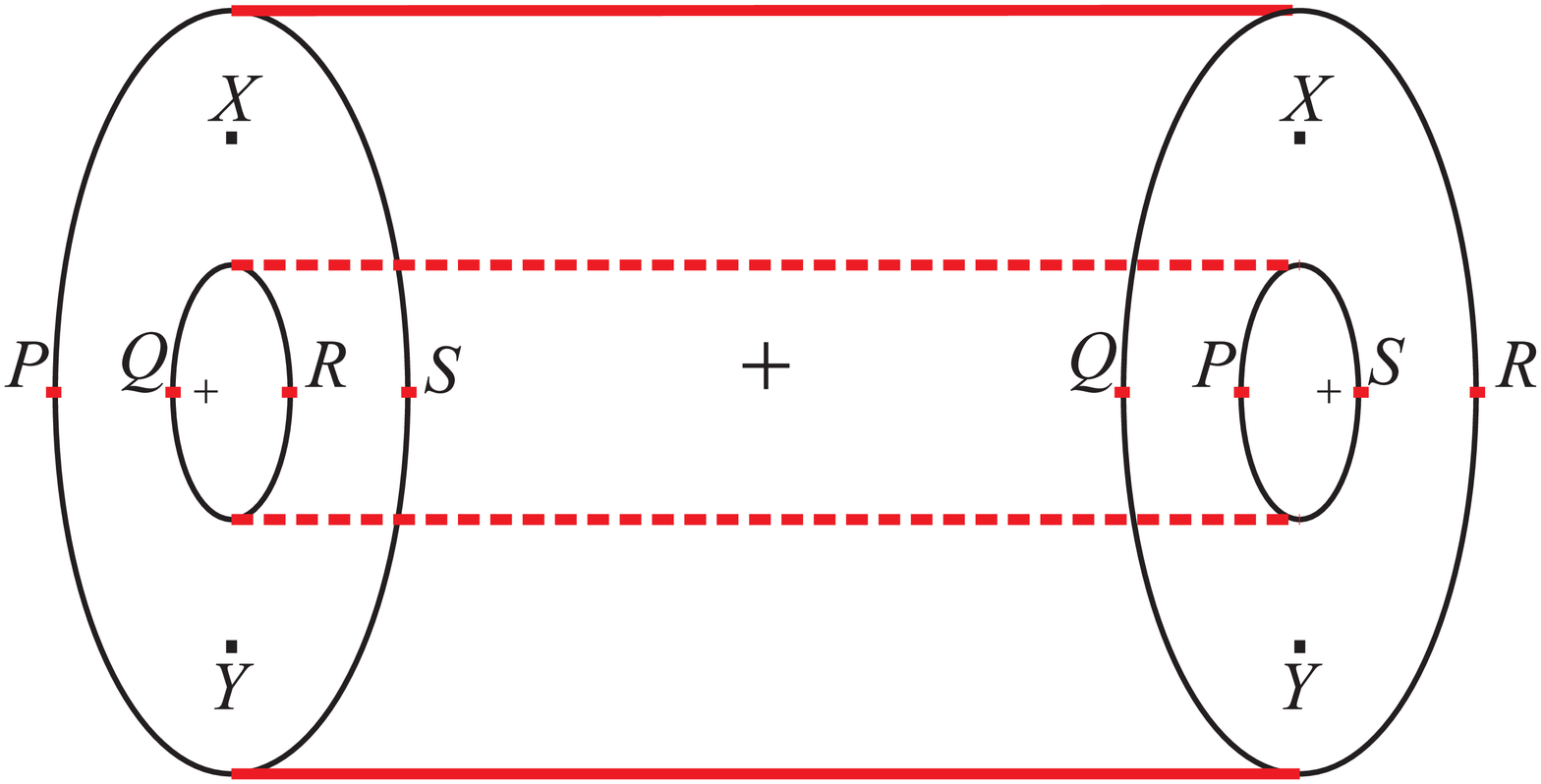}}
\end{center}

\begin{center}
Figure 1.
\end{center}

Let $\xi$ be a tight contact structure on $N$ with minimal convex
boundary of slope $s(T_3)=\infty$. This means that the dividing
curves of $T_3$ consist of two simple closed curves parallel to
$S^1\times \{ pt\}\times\{ 0\}$. Note that the image of $\{
pt\}\times S^1\times I$ under $p$ is an essential annulus in $N$
and the metric closure of its complement in $N$ is a solid torus.
Assume that $T_3$ is a convex torus in standard form with dividing
curves $S^1\times\{0\}\times\{ 0\}$ and
$S^1\times\{\pi\}\times\{0\}$, and $\{ pt\}\times S^1\times\{0\}$,
$pt\in S^1$, are the Legendrian rulings. See Figure 1. The upper
bold line and the upper dashed line form a dividing curve, and the
lower bold line and the lower dashed line form the other dividing
curve. The plus sign $``+"$ in Figure 1 denotes the region
$p(S^1\times [0,\pi]\times \{ 0\})$ in $T_{3}$ bounded by the two
dividing curves.

Let $A$ denote the annulus which is the image of $\{ 0\}\times
S^1\times I$ under $p$. After perturbation, $A$ is convex with
Legendrian boundary. Also assume that $\sharp\Gamma_{A}$, the
number of connected components of the dividing set $\Gamma_A$ of
$A$, is minimal among all convex annuli in its isotopy class
relative to the boundary. $\Gamma_{A}$ contains exactly two
properly embedded arcs. Without loss of generality, we assume that
the endpoints of these two dividing arcs are $P$, $Q$, $R$ and
$S$.

\textbf{Case 5.1.} Both of the two dividing arcs in $\Gamma_{A}$
connect the two different components of $\partial A$.

The two dividing arcs in $\Gamma_A$ must connect the points $P,Q$
and $R,S$ respectively. As shown in Figures 2 and 3, when we cut
$N$ along the convex annulus $A$ and round the edges, we obtain a
solid torus with two dividing curves on the boundary. Moreover,
each of these two dividing curves intersects a meridian of this
solid torus exactly once. There exists a unique tight contact
structure on this solid torus by \cite[Proposition 4.3]{Ho1}. This
implies that in this case, for every choice of $\Gamma_A$, there
exists at most one tight contact structure.

\begin{center}
\scalebox{0.25}[0.25]{\includegraphics {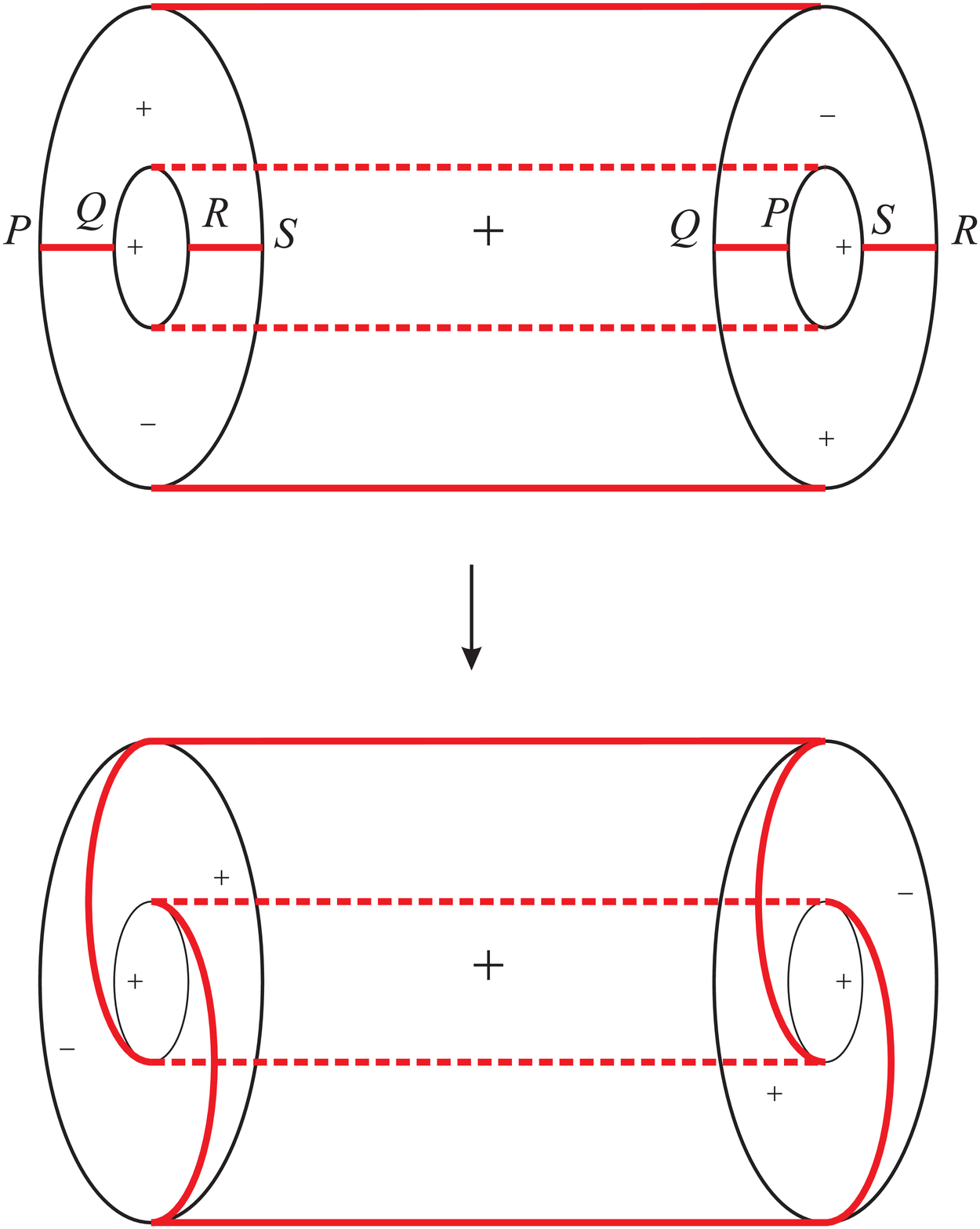}}
\end{center}

\begin{center}
Figure 2.
\end{center}

\begin{center}
\scalebox{0.25}[0.25]{\includegraphics {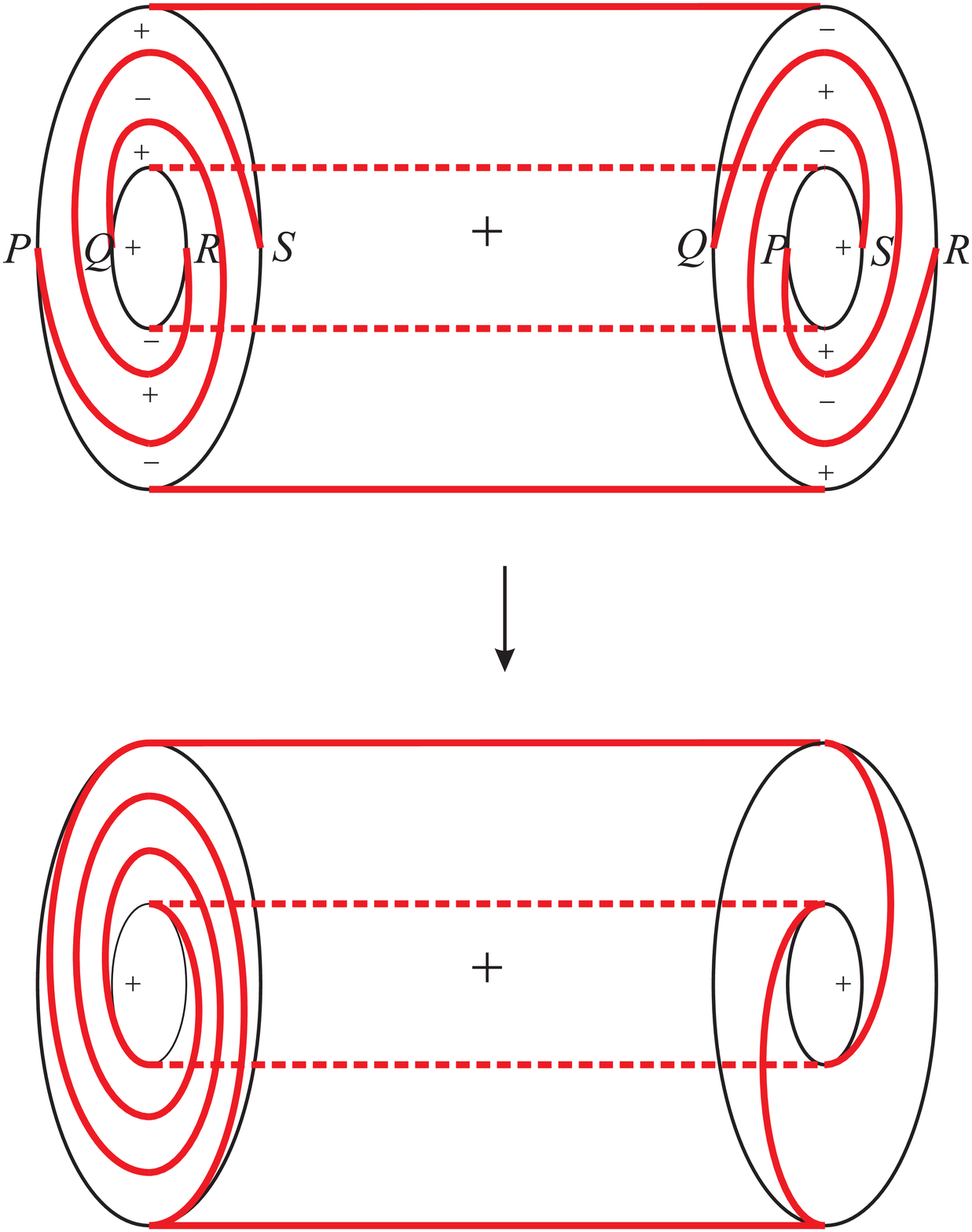}}
\end{center}

\begin{center}
Figure 3.
\end{center}

Similar to the proof of \cite[Proposition 4.9]{Ho1}, we define the
holonomy $k(A)$ by passing to the cover $\{0\}\times
\mathbb{R}\times I\subset S^{1}\times \mathbb{R}\times I$ and
letting $k(A)$ be the integer such that there is a dividing curve
which connects from $(0,\frac{\pi}{2},0)$ to $(0,
2k(A)\pi+\frac{\pi}{2}, 1)$. For example, the holonomy in Figure 2
is $0$, and the holonomy  in Figure 3 is $-1$.

Let $\alpha_0=\cos y dx+\sin y dz$ on $S^1\times S^1\times I$.
Then $\xi_0=\ker \alpha_0$ is the $I$-invariant neighborhood of a
convex $S^1\times S^1$ with dividing curves $S^1\times\{ 0\}$ and
$S^1\times\{ \pi\}$. If we take $\xi_0$ and isotope $S^{1}\times
S^{1}\times \{0\}$ via $(x,y)\longmapsto(x, y-k\pi)$ and isotope
$S^{1}\times S^{1}\times \{1\}$ via $(x,y)\longmapsto(x, y+k\pi)$,
while fixing $S^{1}\times S^{1}\times \{\frac{1}{2}\}$, namely, we
take a self-diffeomorphism of $S^{1}\times S^{1}\times I$ by
sending $(x,y,z)$ to $(x, y+2k\pi(z-\frac{1}{2}), z)$, then we
obtain a tight contact structure $\xi_k$ on $S^{1}\times
S^{1}\times I$ with holonomy $k$ (in the sense of
\cite[Proposition 4.9]{Ho1}), and the corresponding contact form
$\alpha_{k}=\cos
(y+2k\pi(\frac{1}{2}-z))dx+\sin(y+2k\pi(\frac{1}{2}-z))dz$.

Since $\tau^{\ast}(\alpha_{k})=\alpha_{k}$, each nonrotative tight
contact structure $\xi_k$ on $S^{1}\times S^{1}\times I$ is
$\tau$-invariant. So $\xi_{k}$ induces a tight contact structure
on $N$ with holonomy $k(A)=k$. By \cite[Proposition 4.9]{Ho1}, the
nonrotative tight contact structures $\xi_k$, $k\in\mathbb{Z}$, on
$S^{1}\times S^{1}\times I$ are non-isotopic, so they induce
non-isotopic tight contact structures on $N$. All these tight
contact structures on $N$ have Giroux torsion $0$ along $\partial
N$ since each $\xi_k$ has Giroux torsion $0$ along the boundary.
These tight contact structures on $N$ form the subset in Theorem
\ref{number1}(1) whose elements are in 1-1 correspondence with
$\mathbb{Z}$. Note also that for all these tight contact
structures, a convex torus parallel to $\partial N$ must have
slope $\infty$ since each $\xi_k$ is nonrotative.

\textbf{Case 5.2.} The two endpoints of each dividing arc in
$\Gamma_{A}$ belong to the same component of $\partial A$.

If $\Gamma_{A}$ contains an odd number of closed dividing curves,
see Figure 4, then, when we cut $N$ along $A$ and perform
edge-rounding, we find two dividing curves which bound disks. This
contradicts Giroux's criterion (see \cite[Theorem 3.5]{Ho1}). So
$\Gamma_{A}$ must contain an even number of closed dividing
curves.

\begin{center}
\scalebox{0.25}[0.25]{\includegraphics {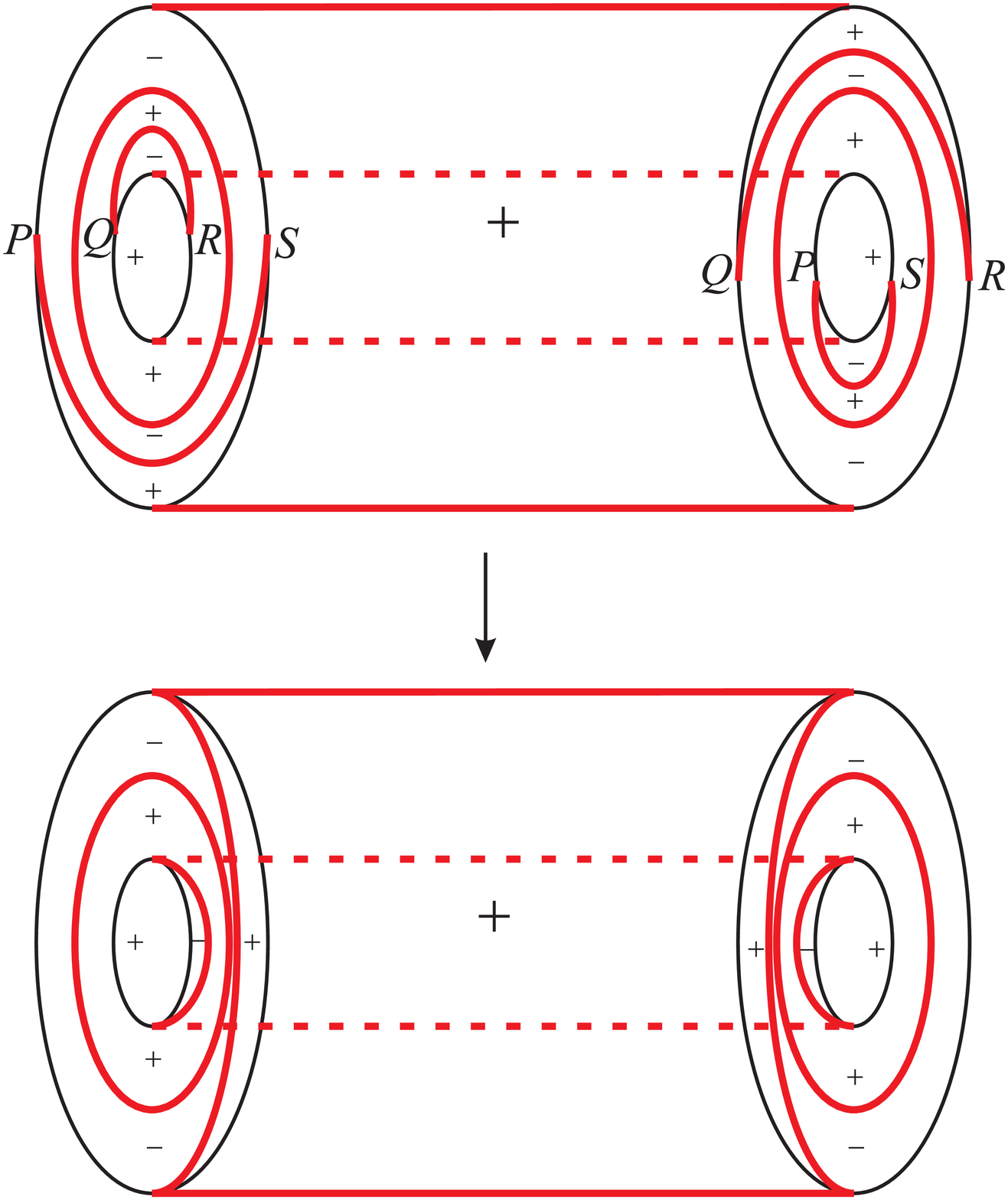}}
\end{center}

\begin{center}
Figure 4.
\end{center}

Suppose $\Gamma_{A}$ contains $2t$ closed dividing curves, where
$t\geq 0$. As shown in Figure 5, when we cut $N$ along the convex
annulus $A$ and round the edges, we obtain a solid torus
$S^1\times D^2$ with $4t+2$ vertical dividing curves.

\begin{center}
\scalebox{0.25}[0.25]{\includegraphics {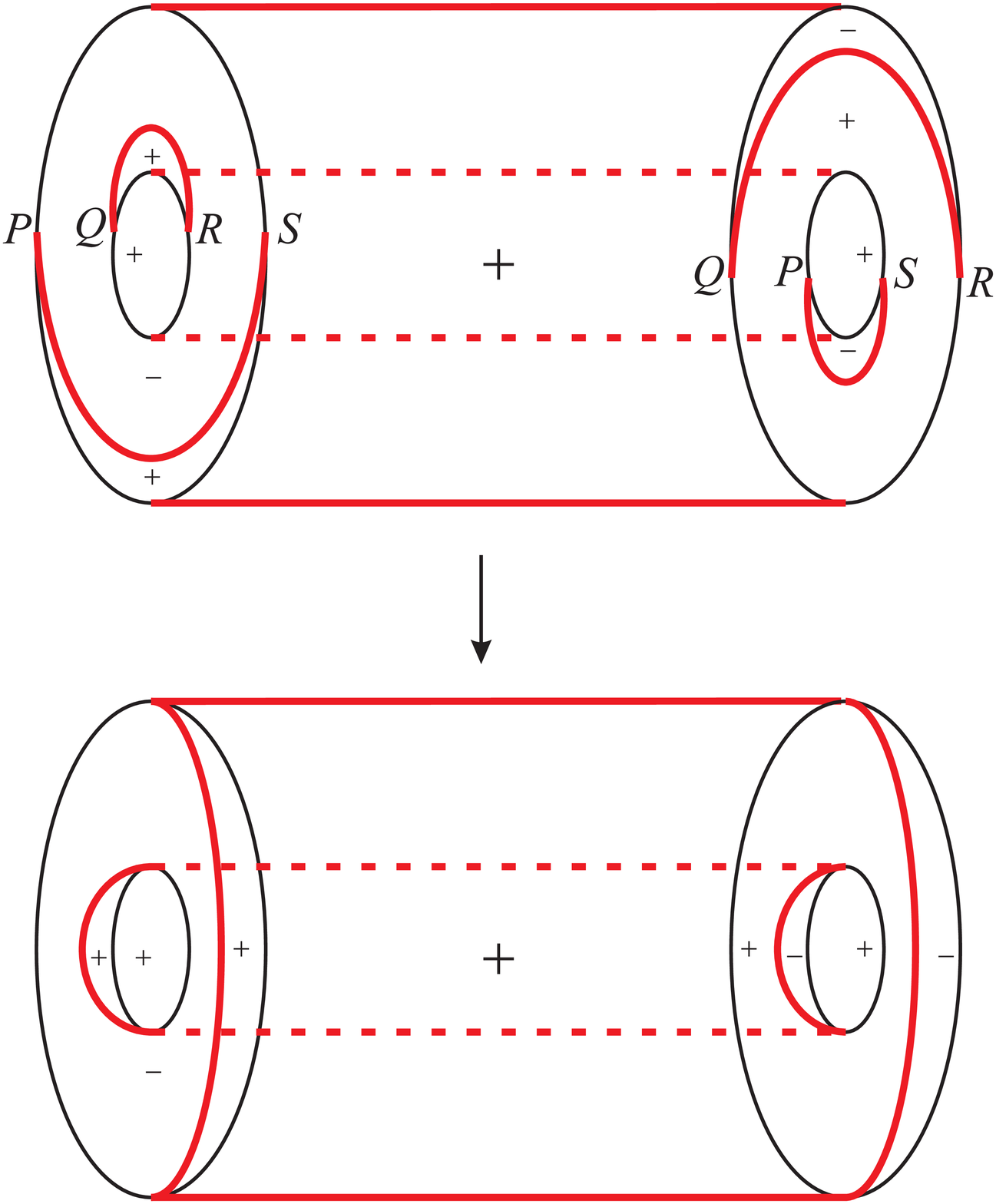}}
\end{center}

\begin{center}
Figure 5.
\end{center}

Next cut $S^1\times D^2$ along a meridional disk $D$ after
modifying the boundary to be standard with horizontal rulings.
Since $\sharp\Gamma_{A}$ is minimal, by a similar argument as in
the proof of \cite[Lemma 5.2]{Ho1}, the dividing set of the convex
meridional disk $D$ has a unique configuration as follows. Let
$\gamma_{0}$ and $\gamma_{1}$ be the two dividing curves on
$\partial (S^1\times D^2)$ which intersect $\partial N$. Then all
$\gamma\in \Gamma_{D}$ must separate $D\cap\gamma_{1}$ from
$D\cap\gamma_{0}$ (hence the dividing curves of $D$ are parallel
segments, with only two boundary-parallel components, each
containing one $D\cap\gamma_i$ in the interior); otherwise there
would exist a bypass which allows for a reduction in the number of
dividing curves on $A$.

Therefore, the tight contact structure $\xi$ on $N$ depends only
on $\Gamma_{A}$, which in turn is determined by the sign of the
boundary-parallel components of $A$ along $\partial N$, together
with $t+2=\sharp\Gamma_{A}$. So in this case, for each $t\ge 0$,
there exist at most two tight contact structures on $N$.

For each $t\in\{0\}\cup\mathbb{Z}^+$, let $\eta_t$ be the contact
structure on $S^{1}\times S^{1}\times I$ given by 1-form
$\beta_{t}=\sin((2t+1)\pi z)dx+\cos((2t+1)\pi z)dy$, with the
boundary adjusted so it becomes convex with two dividing curves on
each component. Let $\eta_t'$ denote the contact structure given
by $-\beta_t$. By \cite[Lemma 5.3]{Ho1}, any two of these tight
contact structures on $S^{1}\times S^{1}\times I$ are distinct.
For each $t\in \{0\}\cup\mathbb{Z}^+$, since
$\tau^{\ast}(\beta_t)=\beta_{t}$, both $\eta_{t}$ and $\eta_{t}'$
are $\tau$-invariant, and hence induce contact structures
$\zeta_t$ and $\zeta_t'$ on $N$ respectively. Since these two
induced contact structures on $N$ lift to distinct tight contact
structures on $S^{1}\times S^{1}\times I$, they are tight and
distinct. Moreover, both $\zeta_t$ and $\zeta_t'$ have minimal
convex boundary of slope $\infty$ and Giroux torsion $t$ along
$\partial N$ by the explicit formula of $\beta_t$ and
\cite[Proposition 3.4]{HKM}.

Similar to \cite[Lemma 5.2]{Ho1} and \cite[Proposition 3.2]{HKM},
if $\Gamma_A$ contains $2t$ closed curves, then $\xi$ is $\zeta_t$
or $\zeta_t'$. $\zeta_0$ and $\zeta_0'$ form the subset in Theorem
\ref{number1}(1) which contains two elements. This completes the
proof of Theorem \ref{number1}(1). For $t\ge 1$, there are exactly
two tight contact structures, namely, $\zeta_t$ and $\zeta_t'$, on
$N$ with minimal convex boundary of slope $\infty$ and Giroux
torsion $t$ along $\partial N$. This proves Theorem
\ref{number1}(2) when $s=\infty$.

Now let $\xi$ be a tight contact structure on $N$ with minimal
convex boundary of slope $s(T_3)=s\in\mathbb{Q}$ and Giroux
torsion $t\geq 1$ along $\partial N$. There is a minimal convex
torus $T'$ in the interior of $N$ which is parallel to $T_{3}$ and
has slope $s$, such that the restriction of $\xi$ on the thickened
torus $U'$ bounded by $T'$ and $T_{3}$ has Giroux torsion $t$.
According to \cite[Lemma 5.2]{Ho1}, $(U', \xi|_{U'})$ is
universally tight.

There is a minimal convex torus $T$ in the interior of $U'$ which
is parallel to $T_{3}$ and has slope $\infty$. We assume that the
restriction of $\xi$ on the thickened torus $U$ bounded by $T$ and
$T_{3}$ is minimally twisting. Note that $U\subset U'$.

The contact submanifold $(N\setminus U, \xi|_{N\setminus U})$
belongs to Case 5.2. If the contact submanifold $(N\setminus U,
\xi|_{N\setminus U})$ belongs to Case 5.1, then each convex torus
in $N\setminus U$ which is parallel to $T$ has slope $\infty$,
contradicting the fact that $T'$ has slope $s$. Note that for the
contact structure $\zeta_0$ on $N$, the slope of a convex torus
parallel to $T_3$ is greater than or equal to $0$. Thus if $s\ge
0$, then the Giroux torsion of $(N\setminus U, \xi|_{N\setminus
U})$ is $t-1$, and if $s<0$, then the Giroux torsion of
$(N\setminus U, \xi|_{N\setminus U})$ is $t$. So there are at most
two possibilities of $(N\setminus U, \xi|_{N\setminus U})$.

Since $U\subset U'$ and $(U', \xi|_{U'})$ is universally tight,
$(U, \xi|_{U})$ is universally tight. By \cite[Proposition
5.1]{Ho1}, there are at most two possibilities of $(U, \xi|_{U})$.
Moreover, these two possibilities are distinguished by their
relative Euler classes. If $(U, \xi|_{U})$ is given, then at most
one possibility of $(N\setminus U, \xi|_{N\setminus U})$ can make
$(N,\xi)$ tight by \cite[Lemma 5.2]{Ho1}. Hence there are at most
two tight contact structures on $N$ with the given conditions.

For a given $t'\in \mathbb{Z}^{+}\cup\{0\}$, let
$0<w<\frac{1}{2t'+1}$ satisfy that $-s=\cot((2t'+1)\pi w)$. Let
$\eta_{t'}$ (we use the same notation as in Case 5.2) be the tight
contact structure on $S^1\times S^1\times [-w,1+w]$ given by
1-form $\beta_{t'}=\sin((2t'+1)\pi z)dx+\cos((2t'+1)\pi z)dy$,
with the boundary adjusted so it becomes convex with two dividing
curves on each component. $\eta_{t'}'$ is given by the 1-form
$-\beta_{t'}$. Think of $N$ as the quotient space of $S^1\times
S^1\times [-w,1+w]$ by identifying $(x,y,z)$ with
$(x+\pi,-y,1-z)$. The transformation $(x,y,z)\mapsto
(x+\pi,-y,1-z)$ on $S^1\times S^1\times [-w,1+w]$ is still denoted
by $\tau$. Since $\beta_s$ is $\tau$-invariant, $\eta_{t'}$ and
$\eta_{t'}'$ induce tight contact structures $\zeta_{t'}$ and
$\zeta_{t'}'$ on $N$ with minimal convex boundary of slope
$s(T_{3})=s$. $\zeta_{t'}$ and $\zeta_{t'}'$ are distinct since
$\eta_{t'}$ and $\eta_{t'}'$ are distinct.

Note that the restriction of the contact structure $\eta_{t'}$ on
$S^{1}\times S^{1}\times [-w,0]$ is minimally twisting. We
conclude that if $s<0$, then the Giroux torsion along $\partial N$
of $\zeta_{t'}$ and $\zeta_{t'}'$ is $t'$ and if $s\geq 0$, then the
Giroux torsion along $\partial N$ of $\zeta_{t'}$ and $\zeta_{t'}'$ is
$t'+1$.

Therefore for each $t\in \mathbb{Z}^+$, there are exactly two
tight contact structures on $N$ with minimal convex boundary of
slope $s$ and Giroux torsion $t$ along $\partial N$. This finishes
the proof of Theorem \ref{number1}(2).

\textbf{Acknowledgements.} The first author is partially supported
by grant no. 10631060 of the National Natural Science Foundation
of China. The second author is partially supported by grant no. 11001171 of
the National Natural Science Foundation of China.

School of Mathematical Sciences, Peking University, Beijing 100871, China

\emph{E-Mail address}:  dingfan@math.pku.edu.cn  \\

Department of Mathematics, Shanghai Jiaotong University, Shanghai
200240, China

\emph{E-mail address}:  liyoulin@sjtu.edu.cn \\

School of Science, Xi'an Jiaotong University,  Xi'an  710049, China

\emph{E-mail address}:  zhangq.math@mail.xjtu.edu.cn

\end{document}